\documentclass[reqno,12pt]{amsart}
\usepackage{amscd,amsfonts,mathrsfs,amsthm,enumerate}
\usepackage{amssymb, amsmath}
\usepackage{stmaryrd}
\usepackage{epsfig}
\usepackage[nobysame]{amsrefs}
\thanks{The author would like to acknowledge support by the General Secretariat for Research and Technology (GSRT) and the
Hellenic Foundation for Research and Innovation (HFRI) Grant No:133.}

\def\real     #1{{\mathbb R^{#1}}}

\def\natural  #1{{\mathbb N^{#1}}}

\def\equationcolor {\color{black}}
\def\textcolor     {\color{black}}

\def\bcoleq    {\begin{equation}\equationcolor}
\def\ecoleq    {\textcolor\end{equation}}
\def\bcoleqn   {\equationcolor\begin{eqnarray}}
\def\ecoleqn   {\end{eqnarray}\textcolor}

\def\gm{{\operatorname{g}_M}}
\def\gn{{\operatorname{g}_N}}

\def\gk{{\operatorname{g}_{M\times N}}}

\def\rind{\operatorname{R}}

\def\dF{\operatorname{d}\hspace{-3pt}F}
\def\df{\operatorname{d}\hspace{-3pt}f}

\def\gind{\operatorname{g}}

\newtheorem{theorem}{Theorem}[section]

\newtheorem*{thm}{Theorem}

\newtheorem{lemma}[theorem]{Lemma}

\theoremstyle{definition}
\newtheorem{remark}[theorem]{Remark}

\newcommand{\bfig}{\begin{figure}}
\newcommand{\efig}{\end{figure}}

\makeatletter
\def\pproof#1{\@ifnextchar[\opargproof
{\opargproof[\it Proof of #1.]}}
\def\opargproof[#1]{\par\noindent {\bf #1 }}

\makeatother

\begin{document}
\title[Minimal maps]{A Schwarz-Pick lemma for minimal maps}
\author[Andreas Savas-Halilaj]{\textsc{Andreas Savas-Halilaj}}
\address{Andreas Savas-Halilaj\newline
University of Ioannina\newline
Section of Algebra \& Geometry \newline
University Campus\newline
45110 Ioannina\newline
Greece\newline
{\sl E-mail address:} {\bf ansavas@uoi.gr}
}
\date{}

\numberwithin{equation}{section}

\begin{abstract}
In this note, we prove a Schwarz-Pick type lemma for minimal maps between negatively curved Riemannian surfaces.
More precisely, we prove that if $f:M\to N$ is a minimal map with bounded Jacobian between two complete negatively curved Riemann surfaces
$M$ and $N$ whose sectional
curvatures $\sigma_M$ and $\sigma_N$ satisfy $\inf\sigma_M\ge\sup\sigma_N$, then $f$ is area decreasing.
\end{abstract}
\maketitle
\section{Introduction}
According to the
Schwarz-Pick Lemma, any non-linear holomorphic map $f:\mathbb{B}\to\mathbb{B}$ from the unit disc $\mathbb{B}$ of the complex plane
$\mathbb{C}$ to itself must be strictly distance decreasing if we endow $\mathbb{B}$ with the Poincar\'e metric.
Ahlfors \cite{ahlfors} exposed in his generalization of the Schwarz-Pick Lemma the
essential role played by the curvature. In the matter of fact, Ahlfors  generalized this lemma to holomorphic mappings between two Riemann
surfaces where the curvature of the domain manifold is bigger than the curvature of the target.
Ahlfors' result was extended by Yau \cite{yau2}.
Yau  showed that if $M$ is a complete K\"ahler manifold with Ricci curvature bounded from below by a constant and $N$ is
another Hermitian manifold with holomorphic bisectional curvature bounded from above by a negative constant, then any
holomorphic mapping from $M$ into $N$ decreases distances up to a constant depending only on the curvatures of $M$ and $N$.
In his proof, Yau exploited a maximum principle at infinity for bounded functions on complete non-compact Riemannian manifolds
with Ricci curvature bounded from below; see for more details \cite{yau1}.

In this paper, we investigate minimal maps between complete Riemann surfaces. According to the terminology introduced by
Schoen \cite{schoen}, a smooth map $f:(M,\gm)\to (N,\gn)$ is called {\it minimal}, if its graph
$$
\Gamma(f):=\big\{(x,f(x))\in M\times N:x\in M\big\}
$$ 
is a minimal surface of the Riemannian product $(M\times N,\gm\times\gn)$. 

There are two important categories of minimal maps between
Riemann surfaces. The first class contains the {\it holomorphic} and {\it anti-holomorphic maps} and the second one the
{\it minimal symplectic maps.}
Eells \cite{eells} proved that a holomorphic or anti-holomorphic map is automatically minimal. Notice that
when both $M$ and $N$ are compact, then there are non-constant holomorphic maps
between them only if the genus of $M$ is greater or equal than the genus of $N$. On the other hand,
the graph of a minimal symplectic map is a minimal Lagrangian
surface. Schoen \cite{schoen} proved an existence and uniqueness result for minimal symplectic diffeomorphisms
between hyperbolic surfaces. i.e., if $\gind_1$ and $\gind_2$ are hyperbolic metrics
on a surface $M$, then there is a unique minimal map $f:(M,\gind_1)\to(M,\gind_2)$ homotopic to the identity map.
If $M$ is a compact hyperbolic surface, due to results of Smoczyk \cite{smoczyk1} and Wang \cite{wang1}
the mean curvature flow deforms a symplectomorphism into a minimal Lagrangian map.

Let us mention here that Aiyama, Akutagawa and Wan \cite{aiyama} obtained a representation formula for a minimal diffeomorphism
between two hyperbolic discs by means of the generalized Gauss map of a complete maximal surface in the anti-de Sitter 3-space.

In this paper, we prove a Schwarz-Pick type lemma for minimal maps between two negatively curved
Riemann surfaces $(M,\gm)$ and $(N,\gn)$. Under natural assumptions, we show that a minimal map from
$M$ to $N$ decreases two dimensional areas. This means that the absolute value $|J_f|$ of the Jacobian determinant
$J_f:=\det(\df)$ of $f$, with respect to the Riemannian metrics $\gm$ and $\gn$, is less or equal than $1$. Such maps are called
{\it area decreasing.} In the case where the Jacobian $J_f$ is identically $1$, the map $f$ is called {\it area preserving.}
\begin{thm}
Let $f:(M,\gm)\to (N,\gn)$ be a minimal map, with bounded Jacobian determinant, between complete negatively curved
Riemannian surfaces
whose sectional curvatures $\sigma_M$ and $\sigma_N$ satisfy
$
\sigma_M\ge-\sigma\ge\sigma_N\ge-\beta,
$
where $\sigma$ is a positive constant. Then $f$ is area decreasing. If, additionally, there exists a point where $f$ is area preserving, then
$\sigma_M=\sigma_N=-\sigma$ and the graph of $f$ is a minimal Lagrangian surface.
\end{thm}

The author and Smoczyk proved in \cite{savas4} that the mean curvature flow of maps between two compact
hyperbolic Riemann surfaces preserves the graphical and the area decreasing property. Due to a result of
Wan \cite{wan} minimal maps between complete Riemann surfaces satisfying the assumptions of our theorem are stable.
Hence, in view of our result, the graphical mean curvature flow can be used to
generate all minimal maps between compact hyperbolic Riemann surfaces.

Finally, let us mention that recently there were proved several Bernstein type results for minimal maps.
In \cite{hasanis1,hasanis2} it is shown that a minimal map $f:\real{2}\to\real{2}$ with bounded
Jacobian determinant must be affine linear. This result was recently generalized by Jost, Xin and Yang in \cite{jost}. 
However, such a result is not true without any assumption on the Jacobian determinant.
For example, the map $f:\real{2}\to\real{2}$ given by
$$
f(x,y)=\frac{1}{2}\big(e^x-3e^{-x}\big)\big(\cos{y}/{2},-\sin{y}/{2}\big)
$$
is a minimal map whose Jacobian determinant takes every value in $\real{}$. Moreover, the graph of $f$ is not holomorphic with
respect to any complex structure of $\real{4}$. Furthermore, due to a result of Torralbo and Urbano \cite{torralbo}, the graph of any minimal
map $f:\mathbb{S}^2\to\mathbb{S}^2$ must be holomorphic or anti-holomorphic.

\section{Maximum principles at infinity}
The root of the maximum principle relies on the following observation: {\it Suppose that $u: M\to\real{}$
is a smooth function
defined on a Riemannian manifold $M$ and assume that it attains at a point $x_0$ a local maximum. Then,}
$$|\nabla u|(x_0)=0\quad\textit{and}\quad\Delta u(x_0)\le 0.$$
Such a point always exists in the case where the manifold is compact. However, the situation is different in the case of complete
non-compact manifolds. To handle the non-compact case,
Omori \cite{omori} proved the following criterion: {\it Suppose
that $M$ is a complete and non-compact Riemannian manifold, with sectional curvatures bounded from below by a constant,
and $u: M\to\real{}$ is a bounded from above smooth function. Then, there exists a sequence of points $\{x_k\}_{k\in\natural{}}$
such that}
\begin{equation}\label{oy}
u(x_k)\ge \sup u-1/k,\quad |\nabla u|(x_k)\le 1/k\quad\textit{and}\quad  \Delta u (x_k)\le 1/k.\tag{O-Y}
\end{equation}
The conclusion \eqref{oy} is known as the {\it Omori-Yau maximum principle}.
Yau \cite{yau1} showed that the conclusion \eqref{oy} holds under the assumption Ricci curvature is bounded from below.

The result of Omori and Yau was generalized by Chen and Xin \cite{chen1} to include cases where the Ricci curvature
may decay at a certain rate. Later Pigola, Rigoli and Setti \cite{pigola} realized that the validity of \eqref{oy}
does not depend on curvature bounds as much as one would expect; for more details see \cite[Theorem 1.9]{pigola}
and \cite{alias}. 
An important case where the Omori-Yau maximum principle can be successfully applied is for properly immersed submanifolds
$F:\Sigma\to L$ with bounded length of the mean curvature vector on a complete Riemannian manifold $L$
with bounded sectional curvatures; see \cite[Example 1.14]{pigola}.

\section{Geometry of graphs}
Let us briefly review some basic facts about the geometry of maps, following closely the presentation
in \cite{savas1,savas4}.
\subsection{Notation}
The graphical submanifold
$$\Gamma(f)=\{(x,f(x))\in M\times N:x\in M\}$$
generated by the smooth map $f:M\to N$ between the Riemann surfaces $(M,\gm)$ and $(N,\gn)$ can be parametrized
via the embedding $F:M\to M\times N$ given by
$F=(I,f),$
where $I:M\to M$ is the identity map. Let us denote by $\pi_M$ and $\pi_N$ the natural projection maps of $M\times N$.
Then, the Riemannian metric $\gk$ on $M\times N$ is given by the formula
$$
\gk=\pi_M^*\gm+\pi_N^*\gn.
$$
We will denote by $\tilde{R}$ the curvature operator of $\gk$. The induced Riemannian metric $\gind$ on the graph $\Gamma(f)$ of $f$ is given by
$$
\gind=\gm+f^*\gn
$$
and its Levi-Civita connection is denoted by $\nabla$. 

Around each point $x\in \Gamma(f)$ we choose an {\it adapted local orthonormal 
frame} $\{e_1,e_2;e_3,e_4\}$
such that $\{e_1,e_2\}$ is tangent and $\{e_3,e_4\}$ is normal to the graph. 
The components of the second fundamental $A$ of the graph with respect to the adapted frame $\{e_1,e_2;e_3,e_4\}$
are denoted as
$$A^{\alpha}_{ij}:=\langle A(e_i,e_j),e_{\alpha}\rangle.$$
Latin indices take values $1$ and $2$ while Greek indices take the values $3$ and $4$. For instance we write
the mean curvature vector in the form
$$H=H^3e_3+H^4e_4. $$
From the {\it Ricci equation} we see that the curvature $\sigma_n$ of the 
normal bundle of $\Gamma(f)$ is given
by the formula
\begin{equation*}
\sigma_n:=\rind^{\perp}_{1234}
=\tilde{R}_{1234}
+A^3_{11}A^4_{12}-A^3_{12}A^{4}_{11}+A^3_{12}A^4_{22}-A^3_{22}A^4_{12}.
\end{equation*}
The sum of the last four terms in the above formula is equal to minus the 
commutator $\sigma^{\perp}$ of the matrices $A^3=(A^3_{ij})$ and $A^4=(A^4_{ij})$, i.e.,
$$
\sigma^{\perp}:=\langle [A^3,A^4]e_1,e_2\rangle=-A^3_{11}A^4_{12}+A^3_{12}A^{4}_{11}-A^3_{12}A^4_{22}+A^3_{22}A^4_{12}.
$$

\subsection{Singular decomposition}\label{singular} There is a natural way to diagonalize the differential $\df$ of $f$. Indeed,
let $\lambda^2\le\mu^2$ be the eigenvalues $f^*\gn$ with respect to $\gm$ at a fixed point $x\in M$. The corresponding
values $0\le\lambda\le\mu$ are called {\it singular values} of $f$ at $x$. Then there exists 
an orthonormal basis $\{\alpha_{1},\alpha_{2}\}$ of
$T_xM$, with respect to $\gm$, and $\{\beta_{1},\beta_{2}\}$ of $T_{f(x)}N$, with respect to the metric
$\gn$, such that
$$\df(\alpha_{1})=\lambda\beta_{1}\quad\text{and}\quad
\df(\alpha_{2})=\mu\beta_{2}.$$
Observe that the vectors
$$v_1:=\frac{\alpha_1}{\sqrt{1+\lambda^2}}\quad\text{and}\quad v_2:=\frac{\alpha_2}{\sqrt{1+\mu^2}}$$
are orthonormal with respect to the metric $\gind$ on the graph of $f$ at $x$. Hence,
\begin{equation*}
e_{1}:=\frac{1}
{\sqrt{1+\lambda^2}}\big(\alpha_1\oplus\lambda\beta_1\big)
\quad\text{and}\quad
e_{2}:=\frac{1}
{\sqrt{1+\mu^2}}\big(\alpha_2\oplus\mu\beta_2\big)\label{tangent}
\end{equation*}
form an orthonormal basis with respect to the metric $\gk$ of the tangent space 
$\dF\left(T_{x}M\right)$ of the graph $\Gamma(f)$ at
$x$. Moreover,
\begin{equation*}
e_{3}:=\frac{1}
{\sqrt{1+\lambda^2}}\big(-\lambda\alpha_1\oplus\beta_1\big)
\quad\text{and}\quad
e_{4}:=\frac{1}
{\sqrt{1+\mu^2}}\big(-\mu\alpha_2\oplus\beta_2\big)\label{normal}
\end{equation*}
form an orthonormal basis with respect to  $\gk$ of the normal space $\mathcal{N}_{x}M$ of the
graph $\Gamma(f)$ at the point $f(x)$.

\subsection{Jacobians of the projection maps}\label{sec 2.4}
Let $\omega_M$ denote the K\"ahler form of the Riemann surface $(M,\gm)$ and $\omega_N$ the K\"ahler form of $(N,\gn)$.
Let us define the parallel forms
$$\omega_1:=\pi^*_M\omega_M\quad\text{and}\quad\omega_2:=\pi^*_N\omega_N.$$
Consider now two smooth functions $u_1$ and $u_2$ given by
$$u_1:=\ast (F^*\omega_1)=\ast\big\{(\pi_M\circ F)^*\omega_M\big\}=\ast (I^*\omega_M) $$
and
$$u_2:=\ast (F^*\omega_2)=\ast\big\{(\pi_N\circ F)^*\omega_N\big\}=\ast (f^*\omega_N) $$
where here $\ast$ stands for the Hodge star operator with respect to the metric $\gind$. Note that $u_1$
is the Jacobian of the projection map from $\Gamma(f)$ to the first factor of $M\times N$ and $u_2$ is the Jacobian
of the projection map of $\Gamma(f)$ to the second factor of $M\times N$.
With respect to the basis $\{e_1,e_2,e_3,e_4\}$ of the singular decomposition, we have
$$u_1=\frac{1}{\sqrt{(1+\lambda^2)(1+\mu^2)}}\quad\text{and}\quad |u_2|=\frac{\lambda\mu}{\sqrt{(1+\lambda^2)(1+\mu^2)}}.$$
The {\it Jacobian determinant} $J_f$ of $f$, with respect to the metrics $\gm$ and $\gn$, is the function defined by the formula
$$f^*\omega_N=J_f\,\omega_M.$$
Therefore,
$$
J_f=\frac{u_2}{u_1}.
$$
The difference $u_1-|u_2|$ measures how far the map $f$ is from being area preserving. In particular,
we say that $f$ is {\it area decreasing} if $u_1-|u_2|\ge 0$ and {\it strictly area decreasing} if $u_1-|u_2|> 0$. If
$u_1-|u_2|\equiv 0$, then $f$ is called {\it area preserving}. Clearly, in the latter situation $f$ is
symplectic.

\subsection{The K{\"a}hler angles}
There are two natural complex structures associated to the product space $(M\times N,\gk)$, i.e.,
$$J_1:=\pi^*_MJ_M-\pi^*_NJ_N\quad\text{and}\quad J_2:=\pi^*_MJ_M+\pi^*_NJ_N,$$
where $J_M$ and $J_N$ are the complex structures on $M$ and $N$ defined by 
$$\omega_M(\cdot\,,\cdot)=\gm(J_M\,\cdot\,,\cdot)\quad\text{and}\quad\omega_N(\cdot\,,\cdot)=\gn(J_N\,\cdot\,,\cdot).$$
Chern and Wolfson in \cite{chern} introduced a function which measures the deviation of
$\dF(T_xM)$ from a complex line of the space $T_{F(x)}(M\times N)$. More precisely, if
we consider $(M\times N,\gk)$ as a complex manifold with respect to $J_1$ then its corresponding
{\it K{\"a}hler angle} $a_1$ is given by
$$\cos a_1=\varphi:=\gk\big(J_1\dF(v_1),\dF(v_2)\big)=u_1-u_2.$$
For our convenience we may require that $a_1\in[0,\pi]$. Observe that, although $\varphi$ is smooth, in general $a_1$
is not smooth at points where $\varphi=\pm 1.$
If there exists a point $x\in M$ where $a_1(x)=0$ then $\dF(T_xM)$ is a complex line of $T_{F(x)}(M\times N)$ and
$x$ is called a {\it complex point} of $F$. If $a_1(x)=\pi$
then $\dF(T_xM)$ is an anti-complex line of $T_{F(x)}(M\times N)$ and $x$ is said {\it anti-complex point} of $F$. In the case
where $a_1(x)=\pi/2$, the point $x$ is called {\it Lagrangian point} of the map $F$. In this case $u_1=u_2$.

Similarly, if we regard $(M\times N,\gk)$ as a K{\"a}hler manifold with respect to the complex structure $J_2$,
then its corresponding K{\"a}hler angle $a_2$ is defined by the formula
$$\cos a_2=\vartheta:=\gk\big(J_2\dF(v_1),\dF(v_2)\big)=u_1+u_2.$$
Notice that $f$ is area decreasing if and
only if both functions $\varphi$ and $\vartheta$ are non-negative. Moreover, observe that there are no points on $M$
where $\varphi=-1$ or $\vartheta=-1$. If $M$ is complete and non-compact, then
$\inf\varphi=-1$ or $\inf\vartheta=-1$ if and only if both singular values $\lambda$ and $\mu$ of $f$ tends to infinity.

\section{Bochner formulas for the Jacobians}
We will derive here the derivative and the Laplacian of a parallel $2$-form on
the product manifold $M\times N$. The proofs are straightforward, make use of the Gauss-Codazzi equations
and can be found in \cite{wang1} (see also
\cite{savas4}). For this reason we omit them.
From now we will always assume that $f$ is a minimal map.

\begin{lemma}\label{fg}
Let $\omega$ be a parallel $2$-form on the product manifold $M\times N$. Then the covariant derivative of the form $F^\ast\omega$ is given by
$$(\nabla_{e_k}F^{\ast}\omega)_{ij}={\sum}_{\alpha}\big(A^{\alpha}_{ki}\omega_{\alpha j}+A^{\alpha}_{kj}\omega_{i\alpha}\big),$$
for any adapted orthonormal frame field $\{e_1,e_2;e_3,e_4\}$.
\end{lemma}

Again by a direct computation we can show the following formula on the Laplacian
of a parallel $2$-form on the product manifold $M\times N$.

\begin{lemma}\label{parallel1}
Let $\Omega$ be a parallel $2$-form on the product manifold $M\times N$. The Laplacian of the form $F^\ast\omega$ is given by the following formula
\begin{eqnarray*}
-(\Delta F^*\omega)_{ij}&=&{\sum}_{\alpha,k,l}\big(A^{\alpha}_{ki}A^{\alpha}_{kl}\omega_{lj}
+A^{\alpha}_{kj}A^{\alpha}_{kl}\omega_{il}\big)
-2{\sum}_{\alpha,\beta, k}A^{\alpha}_{ki}A^{\beta}_{kj}\omega_{\alpha\beta}\\
&&+{\sum}_{\alpha,k}\big(\tilde{R}_{kik\alpha}\omega_{\alpha j}+\tilde{R}_{kjk\alpha}\omega_{i\alpha}\big)
\end{eqnarray*}
where $\{e_1,e_2;e_3,e_4\}$ is an arbitrary adapted local orthonormal frame.
\end{lemma}
From Lemma \ref{parallel1} we can compute the Laplacian of the Jacobians $u_1$ and $u_2$.

\begin{lemma}\label{jacobians}
The Jacobian functions $u_1$ and $u_2$ satisfy the following coupled system of partial differential equations
\begin{eqnarray*}
-\Delta u_1\hspace{-5pt}&=&\hspace{-5pt}\|A\|^2u_1+2\sigma^{\perp} u_2+\sigma_M\big(1-u^2_1-u^2_2\big)u_1-2\sigma_Nu_1u^2_2,\\
-\Delta u_2\hspace{-5pt}&=&\hspace{-5pt}\|A\|^2u_2+2\sigma^{\perp} u_1+\sigma_N\big(1-u^2_1-u^2_2\big)u_2-2\sigma_Mu^2_1u_2.
\end{eqnarray*}
\end{lemma}

Using the special frames introduced in subsection \ref{singular}, from Lemma \ref{fg} and Lemma \ref{jacobians},
by a direct computation we deduce the following:

\begin{lemma}\label{gradlap1}
The gradients of the functions $\varphi$ and $\vartheta$ are given by the equations
$$
2\|\nabla\varphi\|^2=\big(|A|^2-2\sigma^\perp\big)\big(1-\varphi^2\big)\,\,\&\,\,
2\|\nabla\vartheta\|^2=\big(|A|^2+2\sigma^\perp\big)\big(1-\vartheta^2\big).
$$
Moreover, the functions $\varphi$ and $\vartheta$ satisfy the following coupled system of partial differential equations
 \begin{eqnarray*}
-\Delta\varphi&=&\big(|A|^2-2\sigma^{\perp}\big)\varphi+\tfrac{1}
{2}\big(\sigma_M(\varphi+\vartheta)+\sigma_N(\varphi-\vartheta)\big)
(1-\varphi^2),\\
-\Delta\vartheta&=&\big(|A|^2+2\sigma^{\perp}\big)\vartheta+\tfrac{1}{2}\big(\sigma_M(\varphi+\vartheta)-\sigma_N(\varphi-\vartheta)\big)(1-\vartheta^2).
\end{eqnarray*}
\end{lemma}

Observe that away from complex or anti-complex points the second fundamental form quantities
$|A|^2+2\sigma^\perp$ and $|A|^2-2\sigma^\perp$
are expressed  in terms of the cosines of the K\"ahler angles of the graph and of their gradients.

\section{Proof of the theorem}
From our assumptions, the Omori-Yau maximum principle is valid in our setting.
It suffices now  to prove that both $\inf\varphi$ and $\inf\vartheta$ are non-negative numbers. Suppose to the contrary that
$\inf\varphi<0$. Note that since by assumption $J_f$ is bounded, it follows that $\inf\varphi>-1$. Hence from the Omori-Yau
maximum principle we have
that there exists a sequence $\{x_k\}_{k\in\natural{}}$, such that
$$
\lim\varphi(x_k)=\inf\varphi,\quad\lim|\nabla\varphi|(x_k)=0\quad\text{and}\quad\lim\Delta \varphi(x_k)\ge 0.
$$
From Lemma \ref{gradlap1} we have that
\begin{eqnarray*}
-\Delta\varphi(x_k)&=&\frac{2\varphi(x_k)}{1-\varphi^2(x_k)}|\nabla\varphi|^2(x_k)+ \sigma_N(x_k)\varphi(x_k)\big(1-\varphi^2(x_k)\big)\\
&&+\frac{1}{2}\big(\sigma_M(x_k)-\sigma_N(x_k)\big)\big(\varphi(x_k)+\vartheta(x_k)\big)\big(1-\varphi^2(x_k)\big).\end{eqnarray*}
Note that the functions $1-\varphi^2$ and $\varphi+\vartheta$ are positive. Hence, because of
our curvature assumptions the last line of the above equality is non-negative. Passing to the limit we deduce that
$$
0\ge-\sigma\inf\varphi\big(1-(\inf\varphi)^2\big)>0,
$$
which leads to a contradiction. Consequently, $\inf\varphi\ge 0$. Similarly, we prove that $\inf\vartheta\ge 0$.
Hence, the map $f$ must be area decreasing. This completes the first part of the proof.

Let us suppose now that $f$ is an area decreasing map.
Then both $\varphi$ and $\vartheta$ are non-negative functions.
Assume that there is a point $x_0\in M$ where $f$ is area preserving.
Without loss of generality, let assume that $f$ is orientation preserving at $x_0$. Consequently,
$$\varphi(x_0)=0=\min\varphi.$$
From Lemma \ref{gradlap1}, we deduce that
\begin{eqnarray*}
-\Delta\varphi&=&\big\{|A|^2-2\sigma^\perp+\sigma_N(1-\varphi^2)\big\}\varphi\\
&+&\frac{1}{2}\big(\sigma_M-\sigma_N\big)\big(\varphi+\vartheta\big)\big(1-\varphi^2\big)\\
&\ge&2\big\{|A|^2-2\sigma^\perp+\sigma_N(1-\varphi^2)\big\}\varphi.
\end{eqnarray*}
Then from Hopf's strong minimum principle we deduce that $\varphi$ must vanish identically. Going back
to the above identity we obtain that
$$\sigma_M=-\sigma=-\sigma_N,
$$
everywhere. This completes the proof of the theorem.

\begin{remark}
Let us conclude now our paper with some final comments and remarks.
\begin{enumerate}[\rm(a)]
\item
It was very crucial in our proof that the second fundamental terms $|A|^2\pm\sigma^\perp$,
were expressed as gradient terms of the cosines $\varphi$ and $\vartheta$ of the K\"ahler angles of the graph.
However, such a good
structure is not available in higher dimensions and codimensions.
\smallskip
\item
There are various Schwarz-Pick type results for harmonic maps in the literature; see
for instance \cite{chen2,shen,tossati}.
On the other hand, a minimal map $f$ between two Riemannian manifolds $(M,\gm)$ and $(N,\gn)$ becomes harmonic
if we equip $M$ with the graphical metric
$$
\gind=\gm+f^*\gn.
$$
As one can see from the singular value decomposition, the map
$f:(M,\gind)\to (N,\gn)$ is already length decreasing, since its singular values are
$$
\frac{\lambda}{\sqrt{1+\lambda^2}}\quad\text{and}\quad\frac{\mu}{\sqrt{1+\mu^2}}.
$$
Hence, one cannot deduce a Schwarz-Pick type result
for minimal maps by applying directly the already known
results for harmonic maps.
\smallskip
\item
If the map $f$ is holomorphic or anti-holomorphic then, according to the result of Yau \cite{yau2},
it is length decreasing without imposing apriori anything on the size of the differential of $f$.
\end{enumerate}
\end{remark}

\bibliographystyle{y2k}

\end{document}